\documentclass[11pt]{article}

\usepackage{amsmath, amssymb, amsthm, mathtools}
\usepackage{graphicx}
\usepackage{hyperref}
\usepackage[margin=1.05in]{geometry}
\usepackage{enumitem}
\usepackage{booktabs}
\usepackage{listings}
\usepackage{xcolor}      % for custom colours (optional)
\usepackage{longtable}
\usepackage[round]{natbib}
\usepackage{pgfplots}
\pgfplotsset{compat=1.18} % or any recent version

\theoremstyle{definition}

\theoremstyle{remark}

\newcommand{\R}{\mathbb{R}}

\definecolor{codegreen}{rgb}{0,0.6,0}
\definecolor{codegray}{rgb}{0.5,0.5,0.5}
\definecolor{codepurple}{rgb}{0.58,0,0.82}
\definecolor{backcolour}{rgb}{0.95,0.95,0.92}

\begin{document}
%============================================================

\title{Further Improvements to the Lower Bound for an Autoconvolution Inequality}
\author{%
  Aaron Jaech%
  \and
  Alan Joseph\thanks{alanjose@andrew.cmu.edu}%
}
\date{\today}
\maketitle

% ---------------- abstract ----------------
\begin{abstract}

We construct a nonnegative step function comprising 2,399 equally spaced intervals such that  
\[
  \frac{\|f * f\|_{L^{2}(\mathbb{R})}^{2}}{\|f * f\|_{L^{\infty}(\mathbb{R})}\,\|f * f\|_{L^{1}(\mathbb{R})}}
  \;\ge\;
  .926529.
\]
Using a 4× upsampling procedure on this 559-interval optimizer, we further increase the bound to $.94136$,
closing roughly 40\% of the gap between the previous best bound (.901562 on 575 intervals) and the trivial upper limit of 1.

\end{abstract}

\section{Introduction}

A direct application of Hölder’s inequality shows that for any measurable \(F\),
\[
  \frac{\|F\|_{L^{2}(\R)}^{2}}{\|F\|_{L^{\infty}(\R)}\,\|F\|_{L^{1}(\R)}} \le 1,
\]
with equality if and only if \(F\) is an indicator function.  \cite{martin2009supremum} asked how much this bound can be improved when \(F = f * f\) is the autoconvolution of a nonnegative \(f\).  Setting
\[
  c \;=\;\sup_{\substack{f\ge0\\f\in L^1\cap L^2}}
    \frac{\|f*f\|_{L^{2}}^{2}}{\|f*f\|_{L^{\infty}}\,\|f*f\|_{L^{1}}},
\]
one trivially has \(c\le1\).  

Building on the 20-step construction of \cite{matolcsi2010improved}
(\(c\ge.88922\)),
Novikov et al.\ (2025) raised the bound to \(.8962\) using
50 intervals. Very recently, Boyer et al.\ (2025) reported a
575-interval solution (\(c\ge.901562\)) obtained with a
coarse-to-fine gradient ascent.  Independently,
we developed a gradient search that produces a 559-interval function yielding
\[
  \frac{\|f*f\|_{2}^{2}}{\|f*f\|_{\infty}\,\|f*f\|_{1}}
  \;\ge\;
  .926529.
\]
Moreover, by applying a 4× upsampling procedure to our 559-interval solution, we push the lower bound to approximately $c \;\ge\; .94136$,
thereby closing 40\% of the gap to 1 and establishing the new best-known bound on 
$c$.

\section{Numerical Search Algorithm}

Our optimization proceeds in three main stages. We represent a candidate nonnegative step function by its sampled heights \(h\in\R^{N}\) on \(N\) equally spaced intervals, and define the objective
\begin{equation}\label{eq:Cdef}
  C(h)
  \;:=\;
  \frac{\|f_{h} * f_{h}\|_{L^{2}}^{2}}
       {\|f_{h} * f_{h}\|_{L^{1}}
        \,\|f_{h} * f_{h}\|_{L^{\infty}}}\,,
\end{equation}
where \(f\) is the piecewise‐constant function with heights \(h\) supported on \([-1/4,1/4]\).  In practice, we compute \(\|f*f\|_{2}^{2}\) by a Simpson–rule quadrature on the full convolution, \(\|f*f\|_{1}\) by Riemann sum, and \(\|f*f\|_{\infty}\) by its maximum entry. 

\medskip\noindent\textbf{Phase 1: High-LR exploration.}  We initialize a batch of \(B\) random height vectors \(h\sim\mathrm{Unif}[0,1]^{N}\), and optimize them in parallel using the Adam optimizer with learning rate \(3\times10^{-2}\).  At each step \(t\), we add Gaussian gradient noise of standard deviation \(\eta/(t+1)^{\gamma}\) (\(\eta=10^{-3}\), \(\gamma=0.65\)) to promote exploration, then apply a projected update \(h\gets\max(0,h)\).

\medskip\noindent\textbf{Phase 2: Low-LR exploitation.}  Once the exploration phase ends, the learning rate is lowered to \(5\times10^{-3}\) and gradient noise is removed. This fine-tuning phase refines the best candidates found so far.

\medskip\noindent\textbf{Phase 3: Elitist respawn.}  Every \(T\) iterations (e.g.\ \(T=20{,}000\)), we retain the top \(\kappa\%\) of batch members (by current \(C\) value) and replace the remainder with fresh random samples.  This “restart” mechanism prevents stagnation and maintains batch diversity. At the end of the optimization, the best candidate \(h^{*}\) is returned.

\medskip\noindent\textbf{Phase 4: Upsampling and High-Resolution Exploitation.}
Once the best $h^{*}\in\R^{N}$ is found, we apply a simple interpolation‐based upsampling to boost the resolution. We then run gradient‐ascent on this high‐resolution vector, clipping to enforce nonnegativity.

\medskip\noindent\emph{Implementation details.} Our numerical search algorithm is implemented using the JAX library. All experiments were run on an Nvidia A100 GPU via Google Colab, with the full search typically completing in under 10 minutes. We did our search over functions on 768 intervals. After the optimization is complete the leading and trailing zeros were truncated, leaving us with a function on 559 intervals. The source code for our search algorithm is provided in Appendix \ref{app:code}.

\section{Optimized Step Function and Its Autoconvolution}

\begin{figure}[h]
  \centering
  \includegraphics[width=0.8\linewidth]{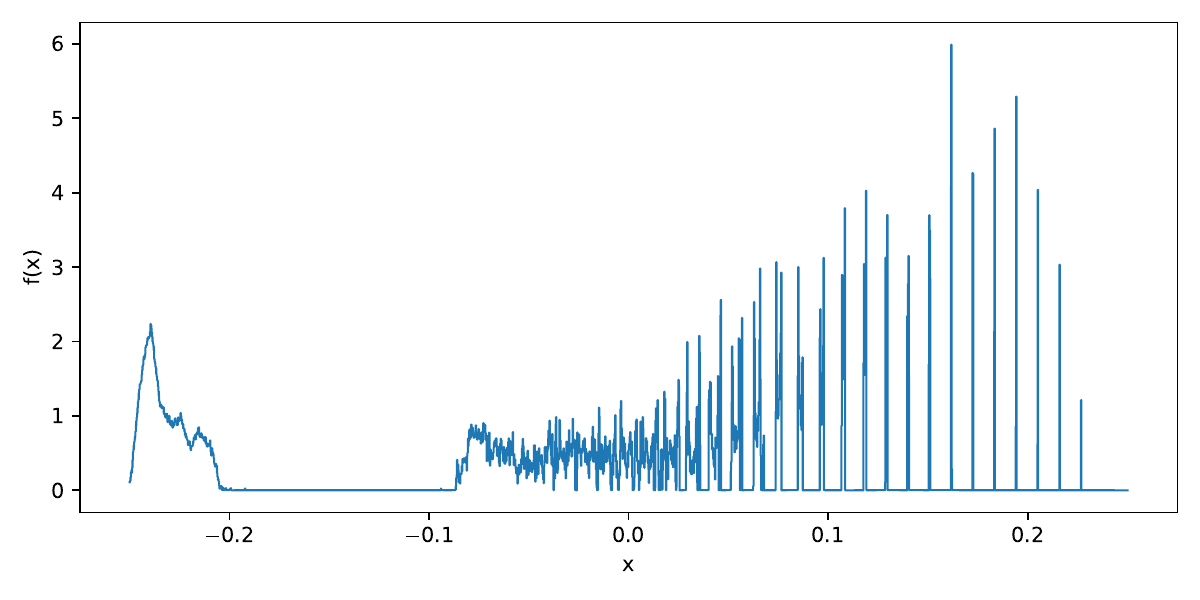}
  \caption{Optimized 559-step function \(f^{*}(x)\).}
  \label{fig:step}
\end{figure}

\begin{figure}[h]
  \centering
  \includegraphics[width=0.8\linewidth]{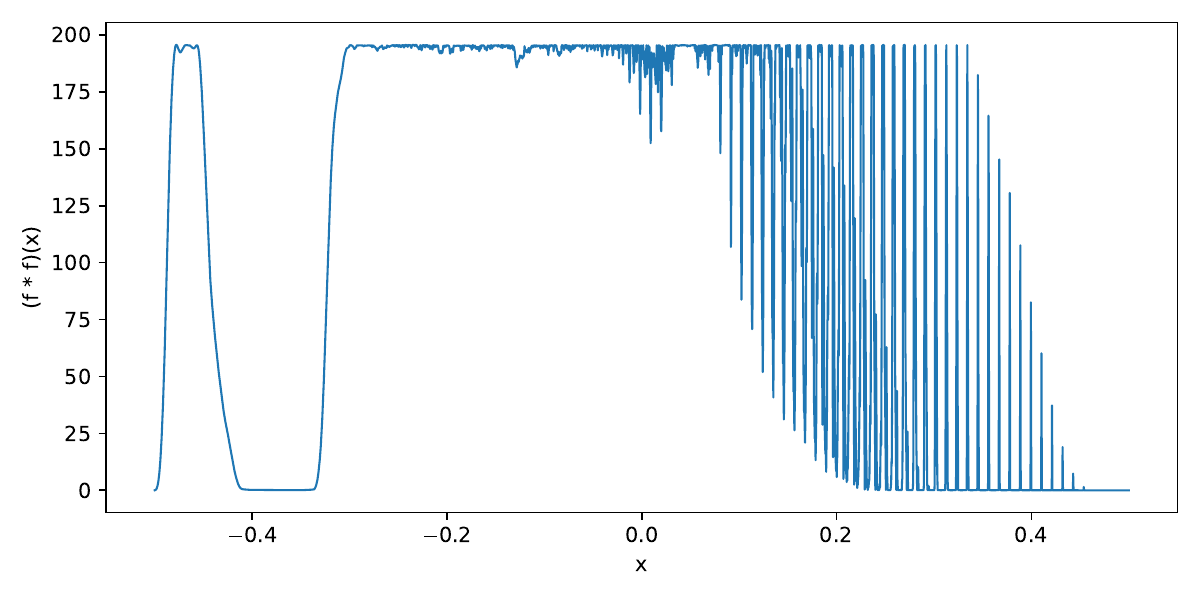}
  \caption{Autoconvolution \(f^{*}*f^{*}(x)\) on \([-1/2,1/2]\).}
  \label{fig:conv}
\end{figure}

Let \(f^{*}\) be the nonnegative 2399-step function achieving $c \ge .94136$.
Figure~\ref{fig:step} plots \(f^{*}(x)\) on \([-1/4,1/4]\): it begins with a very tall, narrow spike near \(x\approx-0.24\), falls sharply to zero, then continues as a “comb” of smaller peaks on the right. The full list of coefficients for the 559-interval function is provided in Appendix \ref{app:heights} and the higher resolution one's coefficients are available on Github.\footnote{\url{https://github.com/ajaech/autocorrelation_inequality}}
Figure~\ref{fig:conv} shows the resulting autoconvolution on \([-1/2,1/2]\): one sees a wide, nearly flat plateau across the center—elevated by \(f^{*}\)’s initial spike—with only gentle ripples appearing on the right edge, directly tracing back to the comb structure.  This arrangement sustains a high baseline to boost the \(L^2\)-norm while keeping both the \(L^1\) and \(L^\infty\) norms under control.

\begin{figure}[h]
  \label{fig:scatter}
  \centering
  \begin{tikzpicture}
    \begin{axis}[
      width=0.7\linewidth,
      xlabel={Number of intervals (log scale)},
      ylabel={Best known lower bound},
      xmode=log,
      log basis x={10},
      xtick={10, 20,50,100,200,500,1000,2000},
      xticklabels={10, 20,50,100,200,500,1000,2000},
      xmin=10, xmax=3500,
      ymin=0.885, ymax=0.95,
    enlarge x limits=false,
         enlarge y limits=false,
      ytick distance=0.01,
      grid=both,
      grid style={dashed,gray!50},
      mark size=2.5pt,
      clip=false                % allows labels to extend past axis box
    ]

    % ---- data points ----
    \addplot+[only marks]
      coordinates {
        (20,  0.88922)
        (50,  0.8962)
        (575, 0.901562)
      };

    % ---- “your” points as red filled triangles ----
    \addplot+[only marks,
      mark=triangle*,
      mark options={fill=red},
      draw=red       % outline color
    ]
      coordinates {
        (50,  0.90331)
        (559, 0.926529)
        (2399,0.9415)
      };

    % ---- labels ----
    \node[font=\footnotesize, anchor=south, xshift=20pt] at (axis cs:20, 0.88922) {Matolcsi--Vinuesa};
    \node[font=\footnotesize, anchor=south] at (axis cs:50, 0.8962)  {Novikov et al.};
    \node[font=\footnotesize, anchor=south] at (axis cs:50, 0.90331) {Our 50-interval};
    \node[font=\footnotesize, anchor=south, xshift=-4pt] at (axis cs:575,0.901562) {Boyer et al.};
    \node[font=\footnotesize, anchor=south,] at (axis cs:559,0.926529) {Our 559-interval};
    \node[font=\footnotesize, anchor=south,  xshift=-18pt] at (axis cs:2399,0.9415) {Our 2399-interval};

    \end{axis}
  \end{tikzpicture}
  \caption{Comparison of published lower bounds versus step-function resolution.}
\end{figure}
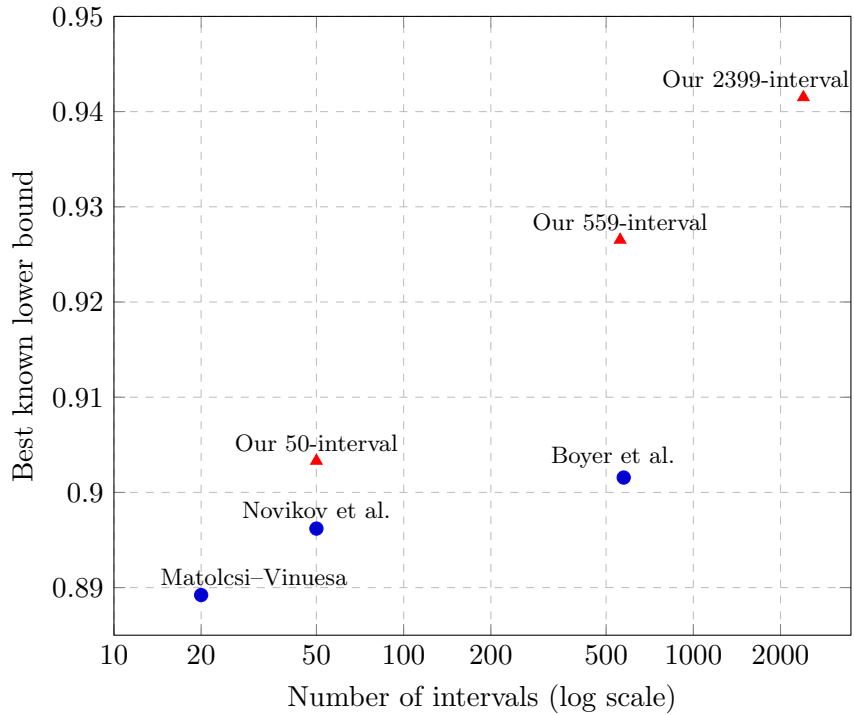

\begin{figure}[h]
  \centering
  \includegraphics[width=0.8\linewidth]{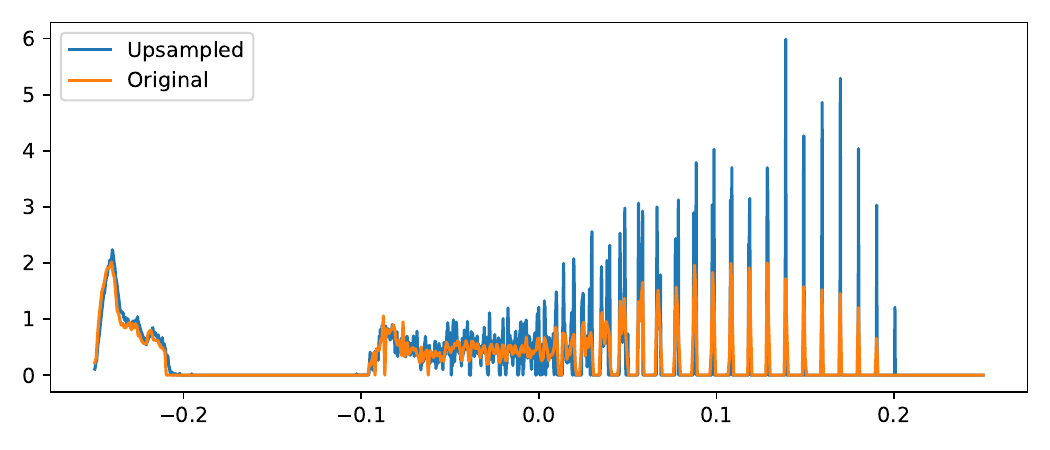}
  \caption{Comparison of the best $f$ before and after optimizing the upsampled version.}
  \label{fig:upsampling}
\end{figure}

\citet{boyer2025improved} were the first to publish a gradient–based numerical search for this inequality.  Their method begins on a coarse grid (23 intervals) and alternates the gradient optimization with upsampling to reach a 575-interval function. Our study was conducted concurrently and follows a complementary strategy: we launch the optimization directly at the target resolution and use Adam’s adaptive learning rate rather than a hand-tuned, decaying step size.  Even when we restrict the search to just 50 intervals, this single-scale approach attains \(.90331\), marginally surpassing the 575-interval value from Boyer and Li. (See Figure \ref{fig:scatter}.)  One plausible explanation is that, once a coarse-grid optimizer has converged, subsequent upsampling may inherit the same local maximum, whereas starting at full resolution (with Adam’s per-coordinate adaptation) lets the search explore fine-scale comb perturbations that appear crucial for narrowing the gap to the upper bound. The comparison of the orange and blue lines in Figure \ref{fig:upsampling} illustrate how the general shape is preserved after gradient ascent is applied to the upsampled function.

\begin{figure}[h]
  \centering
  \includegraphics[width=0.8\linewidth]{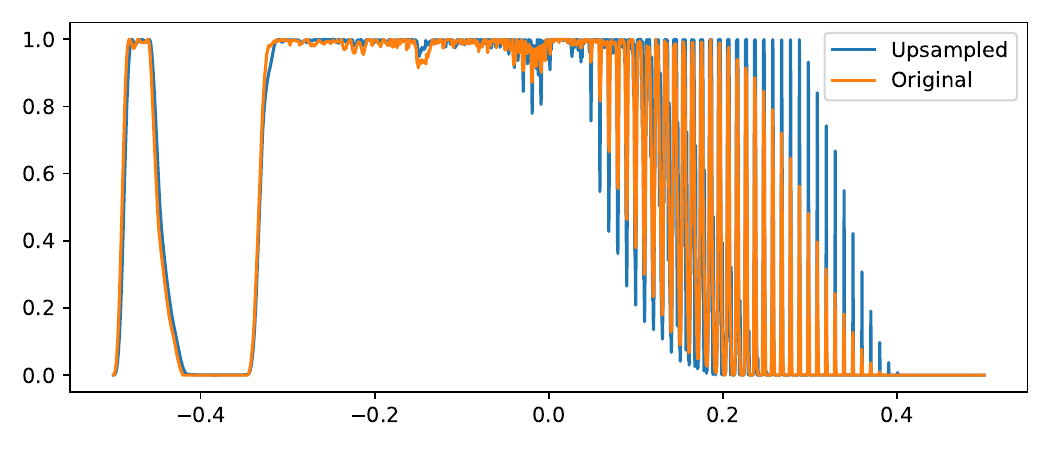}
  \caption{Comparison of $f*f$ before and after optimizing the upsampled version. Both curves are rescaled to have their maximum value fixed at 1.}
  \label{fig:upsampling_conv}
\end{figure}

\section{Conclusion}

We are confident that with longer runs or more exhaustive high-resolution searches the best known lower bound for this inequality can be pushed higher.  
Now that the optimizer consistently discovers a comb-like motif in extremizing step functions, it seems natural to build this structure explicitly into future pipelines, for example through tailored initializations that encourage such fine-scale features from the outset.

Although our study focuses on what \citeauthor{novikov2025alphaevolve} term the \emph{second} autoconvolution inequality, we also applied the same gradient–based search to the \emph{first} and \emph{third} inequalities, which both include a term similar to \(\|f*f\|_{\infty}=\max_{t}(f*f)(t)\) in their numerators.  
%The first,  \(\displaystyle \max_{t}(f*f)(t)\;\ge\;C_{1}\Bigl(\int_{-\tfrac14}^{\tfrac14}f(x)\,dx\Bigr)^{2}\),  and the third,  \(\displaystyle \max_{t}\bigl|(f*f)(t)\bigr|\;\ge\;C_{3}\Bigl(\int_{-\tfrac14}^{\tfrac14}f(x)\,dx\Bigr)^{2}\),  both failed to yield better bounds.  
We attribute this to a \emph{peak‐locking} effect: from the early gradient updates, whichever \(t\) gives the current maximum of \((f*f)(t)\) is immediately reinforced, trapping the search in that local peak.  By contrast, the second inequality’s objective places \(\|f*f\|_{\infty}\) in the denominator, so gradients tend to push down sharp peaks and flatten the convolution globally.  This “peak‐flattening” dynamic prevents premature lock‐in and allows steady improvement toward the best known bound. Devising optimization schemes that navigate such rugged, non-smooth surfaces more effectively remains an open avenue for future work.

\bibliography{main}

\begin{thebibliography}{4}
\expandafter\ifx\csname natexlab\endcsname\relax\def\natexlab#1{#1}\fi

\bibitem[{Boyer and Li(2025)}]{boyer2025improved}
Christopher Boyer and Zane~Kun Li. 2025.
\newblock An improved example for an autoconvolution inequality.
\newblock \emph{arXiv preprint arXiv:2506.16750}.

\bibitem[{Martin and O’Bryant(2009)}]{martin2009supremum}
Greg Martin and Kevin O’Bryant. 2009.
\newblock The supremum of autoconvolutions, with applications to additive number theory.
\newblock \emph{Illinois Journal of Mathematics}, 53(1):219--235.

\bibitem[{Matolcsi and Vinuesa(2010)}]{matolcsi2010improved}
M{\'a}t{\'e} Matolcsi and Carlos Vinuesa. 2010.
\newblock Improved bounds on the supremum of autoconvolutions.
\newblock \emph{Journal of mathematical analysis and applications}, 372(2):439--447.

\bibitem[{Novikov et~al.(2025)Novikov, V{\~u}, Eisenberger, Dupont, Huang, Wagner, Shirobokov, Kozlovskii, Ruiz, Mehrabian et~al.}]{novikov2025alphaevolve}
Alexander Novikov, Ng{\^a}n V{\~u}, Marvin Eisenberger, Emilien Dupont, Po-Sen Huang, Adam~Zsolt Wagner, Sergey Shirobokov, Borislav Kozlovskii, Francisco~JR Ruiz, Abbas Mehrabian, et~al. 2025.
\newblock Alphaevolve: A coding agent for scientific and algorithmic discovery.
\newblock \emph{arXiv preprint arXiv:2506.13131}.

\end{thebibliography}

% ============================================================
\appendix
% ============================================================
\section{559-interval Step function coefficients}\label{app:heights}

\small
\begin{longtable}{llllll}
.22070312 & .30468750 & .72656250 & 1.02343750 & 1.25000000 & 1.48437500 \\
1.57031250 & 1.67968750 & 1.83593750 & 1.90625000 & 1.91406250 & 1.94531250 \\
2.01562500 & 1.82031250 & 1.73437500 & 1.44531250 & 1.14062500 & 1.10156250 \\
.96484375 & .89062500 & .92187500 & .84765625 & .84375000 & .88281250 \\
.92578125 & .93359375 & .80859375 & .91796875 & .83984375 & .91796875 \\
.87500000 & .69531250 & .73437500 & .64062500 & .58984375 & .56250000 \\
.56640625 & .56250000 & .73828125 & .78515625 & .79687500 & .70312500 \\
.62500000 & .62500000 & .61718750 & .62109375 & .54687500 & .49804688 \\
.49023438 & .45898438 & .42382812 & .00000000 & .00000000 & .00000000 \\
.00000000 & .00000000 & .00000000 & .00000000 & .00000000 & .00000000 \\
.00000000 & .00000000 & .00000000 & .00000000 & .00000000 & .00000000 \\
.00000000 & .00000000 & .00000000 & .00000000 & .00000000 & .00000000 \\
.00000000 & .00000000 & .00000000 & .00000000 & .00000000 & .00000000 \\
.00000000 & .00000000 & .00000000 & .00000000 & .00000000 & .00000000 \\
.00000000 & .00000000 & .00000000 & .00000000 & .00000000 & .00000000 \\
.00000000 & .00000000 & .00000000 & .00000000 & .00000000 & .00000000 \\
.00000000 & .00000000 & .00000000 & .00000000 & .00000000 & .00000000 \\
.00000000 & .00000000 & .00000000 & .00000000 & .00000000 & .00000000 \\
.00000000 & .00000000 & .00000000 & .00000000 & .00000000 & .00000000 \\
.00000000 & .00000000 & .00000000 & .00000000 & .00000000 & .00000000 \\
.00000000 & .00000000 & .00000000 & .00000000 & .00000000 & .00000000 \\
.00000000 & .00000000 & .00000000 & .00000000 & .00000000 & .00000000 \\
.00000000 & .00000000 & .00000000 & .00000000 & .00000000 & .00000000 \\
.00000000 & .00000000 & .00000000 & .00000000 & .00000000 & .00000000 \\
.00000000 & .00000000 & .00000000 & .00000000 & .00000000 & .00000000 \\
.00000000 & .00000000 & .00000000 & .00000000 & .00000000 & .00000000 \\
.00000000 & .00000000 & .00000000 & .00000000 & .00000000 & .00000000 \\
.00000000 & .00000000 & .00000000 & .00000000 & .00000000 & .00000000 \\
.00000000 & .00000000 & .00000000 & .00000000 & .00000000 & .00000000 \\
.00000000 & .00000000 & .00000000 & .00000000 & .00000000 & .00000000 \\
.00000000 & .00000000 & .00000000 & .00000000 & .00000000 & .00000000 \\
.00000000 & .00000000 & .00000000 & .00000000 & .00000000 & .00000000 \\
.00000000 & .00000000 & .00000000 & .00000000 & .20898438 & .17773438 \\
.29492188 & .39843750 & .00000000 & .47265625 & .43945312 & .44531250 \\
.61328125 & .67187500 & 1.05468750 & .00000000 & .82812500 & .80078125 \\
.74218750 & .79687500 & .82031250 & .89843750 & .82421875 & .61718750 \\
.53125000 & .62500000 & .51953125 & .34960938 & .94921875 & .33007812 \\
.39453125 & .58203125 & .36914062 & .40429688 & .46875000 & .49023438 \\
.44335938 & .46875000 & .36523438 & .33203125 & .21582031 & .50390625 \\
.24707031 & .35937500 & .40625000 & .45117188 & .00000000 & .43359375 \\
.40625000 & .34570312 & .31445312 & .42968750 & .35156250 & .34570312 \\
.40429688 & .24902344 & .33398438 & .24511719 & .42382812 & .58593750 \\
.52343750 & .46289062 & .40234375 & .49218750 & .47656250 & .56640625 \\
.52734375 & .61328125 & .58203125 & .37109375 & .25390625 & .49023438 \\
.51171875 & .61328125 & .51953125 & .24707031 & .30664062 & .60546875 \\
.59765625 & .48828125 & .47265625 & .56250000 & .49609375 & .30859375 \\
.40820312 & .47656250 & .52734375 & .42187500 & .19433594 & .20800781 \\
.50000000 & .55468750 & .47460938 & .38476562 & .39648438 & .55078125 \\
.45507812 & .21777344 & .31835938 & .53906250 & .51953125 & .25976562 \\
.25781250 & .51171875 & .52734375 & .52734375 & .37109375 & .39453125 \\
.52734375 & .42773438 & .37890625 & .33984375 & .46484375 & .52734375 \\
.40625000 & .42773438 & .67968750 & .32617188 & .48046875 & .31054688 \\
.42968750 & .38085938 & .44140625 & .44335938 & .64453125 & .66015625 \\
.35546875 & .25390625 & .30078125 & .66796875 & .56640625 & .44531250 \\
.25585938 & .36132812 & .34960938 & .39257812 & .51953125 & .85546875 \\
.58593750 & .03662109 & .00000000 & .00000000 & .75000000 & .73828125 \\
.69921875 & .27929688 & .45117188 & .35546875 & .56250000 & .66406250 \\
.73046875 & .14550781 & .00000000 & .00000000 & .00000000 & .36328125 \\
.87500000 & .94140625 & .34179688 & .65234375 & .48632812 & .67968750 \\
.76171875 & .45507812 & .00000000 & .00000000 & .00000000 & .00000000 \\
.00000000 & .75781250 & 1.11718750 & .73437500 & .55468750 & .95703125 \\
.85937500 & .76562500 & .16406250 & .00000000 & .00000000 & .00000000 \\
.00000000 & .00000000 & .16503906 & 1.32031250 & .73828125 & .92187500 \\
1.36718750 & .72265625 & .32226562 & .00000000 & .00000000 & .00000000 \\
.00000000 & .00000000 & .00000000 & .00000000 & 1.31250000 & .96093750 \\
1.35156250 & 1.65625000 & .00000000 & .00000000 & .00000000 & .00000000 \\
.00000000 & .00000000 & .00000000 & .00000000 & .00000000 & 1.09375000 \\
1.51562500 & 1.02343750 & .48437500 & .00000000 & .00000000 & .00000000 \\
.00000000 & .00000000 & .00000000 & .00000000 & .00000000 & .00000000 \\
.52734375 & 1.57031250 & 1.12500000 & .49609375 & .00000000 & .00000000 \\
.00000000 & .00000000 & .00000000 & .00000000 & .00000000 & .00000000 \\
.00000000 & .69921875 & 1.96093750 & 1.48437500 & .00000000 & .00000000 \\
.00000000 & .00000000 & .00000000 & .00000000 & .00000000 & .00000000 \\
.00000000 & .00000000 & .43164062 & 1.82812500 & 1.24218750 & .00000000 \\
.00000000 & .00000000 & .00000000 & .00000000 & .00000000 & .00000000 \\
.00000000 & .00000000 & .00000000 & .33984375 & 1.99218750 & .98437500 \\
.00000000 & .00000000 & .00000000 & .00000000 & .00000000 & .00000000 \\
.00000000 & .00000000 & .00000000 & .00000000 & .00000000 & 1.91406250 \\
.91406250 & .00000000 & .00000000 & .00000000 & .00000000 & .00000000 \\
.00000000 & .00000000 & .00000000 & .00000000 & .00000000 & .00000000 \\
2.00000000 & .84765625 & .00000000 & .00000000 & .00000000 & .00000000 \\
.00000000 & .00000000 & .00000000 & .00000000 & .00000000 & .00000000 \\
.00000000 & 1.71875000 & .71875000 & .00000000 & .00000000 & .00000000 \\
.00000000 & .00000000 & .00000000 & .00000000 & .00000000 & .00000000 \\
.00000000 & .00000000 & 1.57812500 & .30468750 & .00000000 & .00000000 \\
.00000000 & .00000000 & .00000000 & .00000000 & .00000000 & .00000000 \\
.00000000 & .00000000 & .00000000 & 1.52343750 & .00000000 & .00000000 \\
.00000000 & .00000000 & .00000000 & .00000000 & .00000000 & .00000000 \\
.00000000 & .00000000 & .00000000 & .00000000 & 1.45312500 & .00000000 \\
.00000000 & .00000000 & .00000000 & .00000000 & .00000000 & .00000000 \\
.00000000 & .00000000 & .00000000 & .00000000 & .00000000 & 1.20312500 \\
.00000000 & .00000000 & .00000000 & .00000000 & .00000000 & .00000000 \\
.00000000 & .00000000 & .00000000 & .00000000 & .00000000 & .00000000 \\
.65625000 &  &  &  &  &  \\
\end{longtable}

\pagebreak
\section{Source Code for Numerical Search Algorithm}\label{app:code}

\begin{lstlisting}[label={lst:gstep}]
import jax, jax.numpy as jnp, optax, functools, math
from tqdm import trange

def _simpson_integral(y, dx):
    lhs, rhs = y[:-1], y[1:]
    return (dx / 3.0) * jnp.sum(lhs**2 + lhs * rhs + rhs**2)

def _c_single_ref(h):
    h     = jnp.clip(h, 0.0)
    conv  = jnp.convolve(h, h, mode='full')          # length M = 2N-1
    M     = conv.size
    dx    = 1.0 / (M + 1)                            # grid spacing
    # pad with zeros at both ends so we have M+2 samples => M+1 intervals
    y_pad = jnp.concatenate(
        [jnp.zeros(1, conv.dtype), conv, jnp.zeros(1, conv.dtype)])

    l2_sq = _simpson_integral(y_pad, dx)             # ||f*f||_2^2
    l1    = dx * jnp.sum(conv)                       # ||f*f||_1
    linf  = jnp.max(conv)                            # ||f*f||_\infty
    return l2_sq / (l1 * linf)

loss_val_and_grad = jax.vmap(
    jax.value_and_grad(lambda h: -_c_single_ref(jnp.clip(h, 0.0)))
)

explorer_opt, exploiter_opt = optax.adam(3e-2), optax.adam(5e-3)
ETA, GAMMA = 1e-3, 0.55  # noise schedule parameters

@jax.jit
def step_explore(h_raw, opt_state, rng_key, t):
    """High-LR phase with Gaussian gradient noise."""
    vals, grads = loss_val_and_grad(h_raw.astype(jnp.float32))
    sigma       = ETA / ((t + 1) ** GAMMA)
    rng_key, sub = jax.random.split(rng_key)
    grads      += sigma * jax.random.normal(
        sub, grads.shape, grads.dtype)

    updates, new_state = explorer_opt.update(grads.astype(h_raw.dtype),
                                             opt_state, h_raw)
    h_new = optax.apply_updates(h_raw, updates)
    return h_new, new_state, -vals.astype(jnp.float32), rng_key

@jax.jit
def step_exploit(h_raw, opt_state):
    """Low-LR fine-tuning phase (no noise)."""
    vals, grads     = loss_val_and_grad(h_raw.astype(jnp.float32))
    updates, state2 = exploiter_opt.update(grads.astype(h_raw.dtype),
                                           opt_state, h_raw)
    h_new = optax.apply_updates(h_raw, updates)
    return h_new, state2, -vals.astype(jnp.float32)
\end{lstlisting}
\pagebreak
\begin{lstlisting}[title={Gradient-based Search Routine}]
def maximise_c(N=1024, B=2**10, iterations=100000,
               explore_steps=30_000, drop_every=20_000, keep_frac=0.5,
               dtype=jnp.bfloat16, key=jax.random.PRNGKey(0)):
    h_raw   = jax.random.uniform(
                    key, (B, N), minval=0.0, maxval=1.0, dtype=dtype)
    opt_st  = explorer_opt.init(h_raw)

    best_C  = -jnp.inf * jnp.ones(B, jnp.float32)
    best_h  = jnp.clip(h_raw.astype(jnp.float32), 0.0)
    rng     = key

    for t in trange(iterations, desc="optimizing", leave=False):
        if t < explore_steps:
            h_raw, opt_st, C, rng = step_explore(h_raw, opt_st, rng, t)
        else:
            h_raw, opt_st, C      = step_exploit(h_raw, opt_st)

        # update per-candidate bests
        improved = C > best_C
        best_C   = jnp.where(improved, C, best_C)
        best_h   = jnp.where(improved[:, None],
                             jnp.clip(h_raw.astype(jnp.float32), 0.0),
                             best_h)

        # periodic elitist respawn
        if (t + 1) % drop_every == 0:
            K       = int(B * keep_frac)
            top_idx = jnp.argsort(C)[-K:]
            rng, sub = jax.random.split(rng)
            fresh   = jax.random.uniform(
                sub, (B - K, N), minval=0.0, maxval=1.0, dtype=dtype)
            h_raw   = jnp.concatenate([h_raw[top_idx], fresh], axis=0)
            # restart optimizer in the SAME phase we're in
            opt_st  = (
                explorer_opt if t < explore_steps else exploiter_opt
            ).init(h_raw)

    idx = int(jnp.argmax(best_C))
    return best_h[idx], float(best_C[idx]), best_C

best_h, best_C, final_C = maximise_c()
print(f"{best_C:.6f}")
\end{lstlisting}

\pagebreak
\begin{lstlisting}[title={Upsampling Procedure}]
def upsample_1d(h):
    N = h.shape[0]
    x_old = jnp.linspace(-0.5, 0.5, N)
    x_new = jnp.linspace(-0.5, 0.5, 2 * N)
    return jnp.interp(x_new, x_old, h)


def optimize_upsampled(
    best_h, lr=3e-2, iterations=200000, log_every=1000
):
    # ensure float32 for stability
    h = upsample_1d(best_h.astype(jnp.float32))
    # value and gradient of C
    val_and_grad = jax.value_and_grad(
        lambda x: _c_single_ref(jnp.clip(x, 0.0)))
    @jax.jit
    def step(h):
        c, g = val_and_grad(h)
        h_new = jnp.clip(h + lr * g, 0.0)
        return h_new, c
    C_hist = []
    for i in range(iterations):
        h, c = step(h)
        if (i + 1) % log_every == 0:
            C_hist.append(c)
    return h, C_hist

\end{lstlisting}

% ============================================================

\end{document}